\documentclass[a4paper,11pt]{article}
\usepackage{amssymb}
\usepackage{amsmath}
\usepackage{amsthm}
\usepackage{graphicx}
\usepackage{cite}
\usepackage[square,numbers,sort&compress]{natbib}
\newtheorem{thm}{Theorem}[section]
\newtheorem{cor}{Corollary}[section]
\newtheorem{lem}{Lemma}[section]

\theoremstyle{definition}

\theoremstyle{remark}
\newtheorem{rem}{Remark}[section]

\numberwithin{equation}{section}

\begin{document}
	
	\begin{center}
		{\Large \bf Tur\'{a}n type inequalities for rational functions with pescribed poles and restricted zeros}
	\end{center}
	\begin{center}
		{\normalsize \bf Preeti Gupta}
	\end{center}
	\begin{center}
		
			{\normalsize Department of Applied Mathematics\\ Amity University, Noida-201313, India.}\\
			{\normalsize Corresponding auother:  preity8315@gmail.com\\ moblie no. 9717147010 }\\
	
	\end{center}

	\begin{abstract}
	In this paper, we establish some inequalities for rational functions with prescribed poles having s-fold zeros at origin and also show that it implies some inequalities for polynomials and their polar derivatives. 
	\end{abstract}
	
	
\section{Introduction}	
 For each real number $k>0,$ we define the following 
$$D_{k}=\left\lbrace z\in C :\left| z\right|=k  \right\rbrace,$$ 
$$D_{k}^{-}=\left\lbrace z \in C :\left| z\right|<k  \right\rbrace,$$ 
$$D_{k}^{+}=\left\lbrace z\in C:\left| z\right|>k  \right\rbrace.$$ 
If $p(z)$ is a polynomial of degree at most $n$ of a complex variable $z.$ According to the well-known Bernstein's inequality \cite{SB}
\begin{equation}{\label{1.1}}
\max_{z\in D_{1}}| {p}'(z)|\leq n \max_{z\in D_{1}} | p(z) |.	
\end{equation}
The inequality is sharp in the sense that the equality holds if $p(z)=z^{n}.$ Let $P_{n}$ denote the class of all complex polynomials of degree at most $n.$ Let $ \left \| f \right \|=\max_{z\in D_{1}} {|f(z)|},$ the sup-norm of $f$ on the unit circle. Then Bernstein's inequality can be restated as the following extremal problem:
$$\max_{p\in P_{n}} {\frac{||p^{'}||}{||p||}}=n.$$

If we restrict ourselves to the class of polynomials having all its zeros in $D_{1}\cup D_{1}^{+},$ then it was conjectured by Erd\"{o}s and later verified by Lax \cite{PD}
\begin{equation}{\label{1.2}}
|| {p}'||\leq \frac{n}{2}||p||.		
\end{equation}
On the other hand, if all the zeros of $p\in P_{n} $ lie in $D_{1}\cup D_{1}^{-},$ then it was proved by Tur$\acute{a}$n \cite{PT} that
\begin{equation}{\label{1.3}}
||{p}'||\geq \frac{n}{2}||p||.	
\end{equation}
In the literature \cite{GTM,QG,PB}, there exist several improvement and generalization of inequalities (\ref{1.1}), (\ref{1.2}) and (\ref{1.3}).\\
\section{Rational functions}
Let $a_{1},a_{2},\dots,a_{n} $ be $n$ given points in $|z|>1.$ We will consider the following space of rational functions with prescribed poles:
$$\Re_{n}=\Re_{n}\left(a_{1},a_{2},\dots,a_{n} \right):= \left\lbrace \frac{p(z)}{w(z)} : p \in P_{n} \right\rbrace ,$$
where $w(z):=(z-a_{1})(z-a_{2})\dots(z-a_{n}). $\\

Denote
$$B(z)=\frac{w^{*}(z)}{w(z)}=\frac{z^{n}\overline{w\left ( \frac{1}{\bar{z}} \right )}}{w(z)}
:=\prod_{j=1}^{n}\left(\frac{1-\bar{a}_{j}z}{z-a_{j}} \right).$$
The product $B(z)$ is known as Blaschke product and one can easily show that $ \left|B(z) \right| =1 $ and $\frac{zB^{'}(z)}{B(z)}=|B^{'}(z)| $ for $ z\in D_{1} .$\\

The inequalities of Bernstein and  Tur$\acute{a}$n have been extended to the rational functions Li, Mohapatra and Rodriguez \cite{XRR}:\\
                 If $|z|=1,$ then for any $r\in \Re_{n},$
	\begin{align}{\label{2.1}}
\left|{r}'(z) \right| \leq\left|{B}'(z) \right|\left\|r\right\|,		
	\end{align} 

	\begin{align}{\label{2.2}}
		\left|{r}'(z) \right| \geq \frac{1}{2}\left|{B}'(z) \right|\left\|r \right\|.	
	\end{align}
  Aziz and Shah \cite{AW97,AW} proved the following theorems which improve upon the inequality (\ref{2.2}).
 
  	\begin{thm}{\label{thm2.1}}
  	Let $ r\in \Re_{n}, $ where r has exactly n poles at $ a_{1},a_{2},\dots,a_{n} $ and all its zeros lie in $D_{1}\cup D_{1}^{-}.$ Then for $z\in D_{1}$
  	\begin{align}{\label{2.3}}
  		\left|{r}'(z) \right| \geq\frac{1}{2}\left|{B}'(z) \right|\left[ \left|r\left(z \right) \right| +m^{'} \right] , 
  	\end{align} 
  	where $m^{'}=\min_{z\in D_{1}}\left|r(z) \right|.$ 
  	Equality attains for $r(z)=B(z)+he^{i\alpha }$ with $h\leq1$ and $\alpha$ is real.	
  \end{thm}
\begin{thm}{\label{thm2.2}}
	If $ r\in \Re_{n} $ has all its zeros in $D_{k}\cup D_{k}^{-},$ $k \leq 1, $ then for $z \in D_{1},$ we have 
	\begin{align}{\label{2.4}}
		\left|{r}'(z) \right| \geq\frac{1}{2}\left\lbrace \left|{B}'(z)\right|-\frac{n(1+k)-2m}{1+k} \right\rbrace \left|r(z) \right|  , 
	\end{align} 
	where m is the number of zeros of $r.$ 
\end{thm}
As a generalization of Theorem \ref{thm2.2}, Wali \cite{wali} obtained the following result.

\begin{thm}{\label{thm2.3}}
	Suppose $ r\in \Re_{n}, $ where r has exactly n poles at $ a_{1},a_{2},\dots,a_{n} $ and all its zeros lie in $D_{k}\cup D_{k}^{-},$ $k \leq 1,$ with a zero of order $s$ at origin.Then for $z\in D_{1}$
	\begin{align}{\label{2.5}}
		\left|{r}'(z) \right| \geq\frac{1}{2}\left\lbrace \left|{B}'(z)\right|- n + \frac{2(m+sk)}{1+k} \right\rbrace \left|r(z) \right|  ,  
	\end{align} 
	where m is the number of zeros of $r.$
	
\end{thm}

 \section{Main results}	
  We present the following results, which provides generalizations as well as refinements of above inequalities (\ref{2.3}) and  (\ref{2.5}) . Infact, we prove the following:

\begin{thm}{\label{thm3.1}}
	Suppose $r(z)=\frac{p(z)}{w(z)} \in \Re_n,$ where $r$ has exactly $n$ poles at $ a_{1},a_{2},\dots,a_{n} $ and all its zeros lie in $D_{k}\cup D_{k}^{-},$ $k \leq 1,$ with a zero of order $s$ at origin, that is, $p(z)=c_{m}\prod_{j=1}^{m-s}(z-b_{j})=c_{0}+c_{1}z+\dots+c_{m}z^{m-s},$ $c_{m} \neq 0 ,$ $ m \leq n ,$ $ |b_{j}|\leq 1,$ $  k \leq 1 , j=1,2,\dots, m-s.$  Then for $z\in D_{1}$
	\begin{align}{\label{3.1}}
	\left|{r}'(z) \right| \geq\frac{1}{2}\left\lbrace \left|{B}'(z)\right|- n + \frac{2(m+sk)}{1+k} +2\left(\sum_{j=1}^{m-s}\frac{1}{1+|b_{j}|}-\frac{m-s}{1+k} \right)  \right\rbrace \left|r(z) \right|  , 
\end{align} 
where m is the number of zeros of $r.$ The result is sharp and equality holds for $r(z)=\frac{z^{s}(z+k)^{m-s}}{(z-a)^{n}}$ and $B(z)=(\frac{1-az}{z-a})^{n}$ at $z=1$ and $m \leq n,$ $k \leq 1 < a.$

\end{thm}
\begin{rem}
	Since under the condition of Theorem \ref{thm3.1} 
	$$\sum_{j=1}^{m-s}\frac{1}{1+|b_{j}|}-\frac{m-s}{1+k} \geq 0, $$
 One can see at once that Theorem \ref{thm3.1} improves  inequality (\ref{2.5}). 
\end{rem}

For $s=0,$ Theorem \ref{thm3.1} reduces to Rather, Iqbal and Dar \cite{NAI}.\\
For $m=n,$ Theorem \ref{thm3.1} reduces to the following result.
\begin{cor}
Let $ r\in \Re_n, $ where r has exactly n poles at $ a_{1},a_{2},\dots,a_{n} $ and all its zeros lie in $D_{k}\cup D_{k}^{-}$, $k \leq 1,$  with a zero of order $s$ at origin then for $z\in D_{1},$
\begin{align}{\label{3.2}}
	\left|{r}'(z) \right| \geq \frac{1}{2}\left[\left|{B}'(z) \right|+ \frac{n(1-k)+2sk}{1+k} +2\left( \sum_{j=1}^{n-s}\frac{1}{1+|b_{j}|}-\frac{n-s}{1+k}  \right)  \right] \left|r(z) \right|.	
\end{align}	
\end{cor}
\begin{cor}{\label{cor3.2}}
	Suppose $r(z)=\frac{p(z)}{w(z)} \in \Re_n,$ where $r$ has exactly $n$ poles at $ a_{1},a_{2},\dots,a_{n} $ and all its zeros lie in $D_{k}\cup D_{k}^{-},$ $k \leq 1,$ with a zero of order $s$ at origin, that is, $p(z)=c_{m}\prod_{j=1}^{m-s}(z-b_{j})=c_{0}+c_{1}z+\dots+c_{m}z^{m-s},$ $c_{m} \neq 0 ,$ $ m \leq n ,$ $ |b_{j}|\leq 1,$ $  k \leq 1 , j=1,2,\dots, m-s.$  Then for $z\in D_{1}$
	\begin{align}{\label{3.3}}
		\left|{r}'(z) \right| \geq\frac{1}{2}\left\lbrace \left|{B}'(z)\right|- n + \frac{2(m+sk)}{1+k}+ \frac{2k}{k+1}\left( \frac{k^{m}|c_{m}|-|c_{s}|}{k^{m}|c_{m}|+|c_{s}|}\right)   \right\rbrace \left|r(z) \right|  , 
	\end{align} 
	where m is the number of zeros of $r.$ The result is sharp and equality holds for $r(z)=\frac{z^{s}(z+k)^{m-s}}{(z-a)^{n}}$ and $B(z)=(\frac{1-az}{z-a})^{n}$ at $z=1$ and $m \leq n,$ $k \leq 1 < a.$
	
\end{cor}

In Corollary \ref{cor3.2} if we take  $k=1,$ we get Theorem 2.1 of Wali \cite{wali}. 
\begin{rem}
	Since under the hypothesis of Corollary \ref{3.2} $$k^{m}|c_{m}|-|c_{s}| \geq 0 $$  We can see that the bound of (\ref{3.3}) is an improvement of Theorem \ref{thm2.3}.
	\end{rem}
If we take $s=0$ in Corollary \ref{cor3.2}, we get the following result.
\begin{cor}{\label{cor3.3}}
Suppose $r \in \Re_n,$ where $r$ has $n$ poles and all the zeros of $r$ lie in  $D_{k}\cup D_{k}^{-},$ $k \leq 1.$ Then for $z\in D_{1}$
	\begin{align}{\label{3.5}}
	\left|{r}'(z) \right| \geq\frac{1}{2}\left\lbrace \left|{B}'(z)\right|- n + \frac{2m}{1+k}+ \frac{2k}{k+1}\left( \frac{k^{m}|c_{m}|-|c_{0}|}{k^{m}|c_{m}|+|c_{0}|}\right)   \right\rbrace \left|r(z) \right|  
\end{align} 
\end{cor}
If we put $k=1,$ in  inequality (\ref{3.5}), we get Corollary 2.2 of Wali \cite{wali}. \\
In case number of poles of $r$ is same as its zeros, that is, when $m=n$ then Corollary \ref {cor3.3} gives an improvement of inequality of Aziz and Shah \cite{AW}.\\

We also discuss some consequences of Theorem \ref{thm3.1}. If we consider $p(z)$ as a polynomial of degree $n.$ Let us take $a_{j}=\alpha,$ $ |\alpha|\geq 1,$ for $j=1,2,\dots,n,$ then $w(z)=(z- \alpha)^{n}$ and $r\left(z \right)=\frac{p\left(z \right) }{\left(z- \alpha \right)^{n} },$ and hence we have		
\begin{align}{\label{3.6}}
	{r}'\left(z \right)
	&=\frac{(z- \alpha)^{n}{p}'(z)-n(z-\alpha)^{n-1}p(z)}{(z-\alpha)^{2n}}\nonumber\\	
	&=-\left\lbrace \frac{np(z)-(z- \alpha){p}'(z)}{(z-\alpha)^{n+1}} \right\rbrace =\frac{-D_{\alpha}p(z)}{(z- \alpha)^{n+1}},
\end{align}
where $D_{\alpha}p(z) = np(z) + (\alpha - z ){p}'(z) $ is the polar derivative of $p(z)$ with respect to point $\alpha.$ It generalizes the ordinary derivative ${p}'(z)$ of $p(z)$ in the sense that 
$$ \lim_{\alpha \rightarrow\infty}\frac{D_{\alpha}p(z)}{\alpha}={p}'(z).$$
For $B\left(z \right)=\frac{w^{*}\left(z \right) }{w\left(z \right) }, $ we have 
${B}'\left(z \right)=\frac{n\left(\left|\alpha \right|^{2} -1 \right) }{\left( z- \alpha\right) ^{2}}\left(\frac{1-\bar{\alpha}z}{z- \alpha} \right) ^{n-1}.$ Hence
for $z \in T_{1},$   $\left|{B}'\left(z \right) \right|=\frac{n\left(\left|\alpha \right| ^{2}-1 \right) }{\left|z- \alpha \right| ^{2}}.$\\

With these choices, and \label{key} from  Corollary \ref{cor3.2} for $|z|=1,$ we obtain the following result in terms of polar derivative.
By letting $\left|\alpha \right| \rightarrow \infty,$ we have the following result which generalizes Tur$\acute{a}$n's inequality \cite{PT}.
\begin{cor}{\label{cor3.4}}
	Suppose $p(z)$ is a polynomial of degree $n,$ having  all its zeros lie in $D_{k}\cup D_{k}^{-},$ $k \leq 1,$ with a zero of order $s$ at origin, that is, $p(z)=z^{s}h(z)$ where $h(z)=c_{n}\prod_{j=1}^{n-s}(z-b_{j})=c_{0}+c_{1}z+\dots+c_{n}z^{n-s},$ $c_{n} \neq 0 , |b_{j}|\leq 1,  k \leq 1 ,j=1,2,\dots, n-s.$  Then for $z\in D_{1}$
	\begin{align}{\label{3.7}}
		\left|{p}'(z) \right| \geq\frac{n}{k+1}\left\lbrace 1+ \frac{sk}{n}+ \frac{k}{n}\left( \frac{k^{n}|c_{n}|-|c_{s}|}{k^{n}|c_{n}|+|c_{s}|}\right)   \right\rbrace \left|p(z) \right| .
	\end{align} 
	
\end{cor}

\begin{rem}
	The inequality (\ref{3.7}), is refines a well known polynomial inequality due to Malik \cite{malik}.
\end{rem}

	\begin{thm}{\label{thm3.2}}
	Let $r\left( z \right)= \frac{p\left( z\right) }{w\left( z\right) } \in \Re_{n}$ and $b_{1},b_{2},\dots,b_{m}$ are be zeros of $r\left( z \right)$ all lying in $D_{k}\cup D_{k}^{-},$ $k\leq 1$. Then for $z \in D_{1},$ 
	$$ 	\left|{r}'(z) \right|\geq \frac{1}{2}\left[\left|{B}'\left(z \right)  \right|  + \frac{n(1-k)-2m}{1+k}+ 2\left(\sum_{j=1}^{m}\frac{1}{1+\left|b_{j} \right| }-\frac{n-m}{1+k} \right)     \right]\left(\left|r(z) \right| +m^{'} \right), $$ 
	where $m$ is the number of zeros of $r$ and $m^{'}=\min_{z\in T_{k}}\left|r(z) \right|$.
\end{thm}
If $r\left( z \right) $ has exactly $n$ zeros all lying in $D_{k}\cup D_{k}^{-},$ where $k \leq 1,$  we get Theorem 2.2 of Hans and Mir \cite{SA}.

	\section{Lemmas}	
For the proof of theorems we require the following lemmas. We begin by starting the following lemma is due to Li , Mohapatra and Rodriguez ~\cite{XRR}.
\begin{lem}{\label{4.1}}
	If $r \in \Re_{n}$ and $r^{*}\left ( z \right )= B\left ( z \right )\overline{r\left ( \frac{1}{\bar{z}} \right )}$ then for $z\in T_{1}$, we have 
	\begin{align*}
		\left|{\left(r^{*}(z) \right)}'  \right|+\left|{r}'(z) \right| \leq \left|{B}'(z) \right|\left \| r \right \|.   
	\end{align*}
\end{lem}
Next lemma is due to Bidkham and Shahmansouri \cite{MT}.
\begin{lem}{~\label{4.2}}
If $z \in T_{1}$, then 
$$Re\left( \frac{z{w}'(z)}{w(z)}\right)=\frac{n-\left|{B}'(z) \right| }{2}. $$	

\end{lem}
\begin{lem}
If  be a sequence of real numbers such that $0 \leq x_{j} \leq 1 ,$ $j \in N$ then
$$\sum_{j=1}^{n}\frac{1-x_{j}}{1+x_{j}} \geq \frac{1-\prod_{j=1}^{n}x_{j}}{1+\prod_{j=1}^{n}x_{j}}, \  \text{for all} \  n \in N$$
\end{lem}
A direct consequence of the above lemma we get the following result.
\begin{lem}{~\label{4.4}}
	If  be a sequence of real numbers such that $ x_{j} \geq 1, $ $j \in N$ then
	$$\sum_{j=1}^{n}\frac{1-x_{j}}{1+x_{j}} \leq \frac{1-\prod_{j=1}^{n}x_{j}}{1+\prod_{j=1}^{n}x_{j}}, \  \text{for all} \  n \in N$$
\end{lem}

\section{Proof of theorems}	
\begin{proof}[\textbf{Proof of Theorem 3.1}] For $z \in D_{1}$ and $r(z)=\frac{z^{s}p(z)}{w(z)}$ with  $p(z)=c_{m}\prod_{j=1}^{m-s}(z-b_{j}),$ $c_{m-s}\neq 0$, $m\leq n,$ $|b_{j}|\leq k\leq 1,$ $ j=1,2,\dots,m-s.$
	\begin{align*}
	  Re\left( \frac{z{r}'\left(z \right) }{r\left(z \right) }\right) 
		&=	s+Re\left(  \frac{z{p}'\left(z \right) }{p\left(z \right) } \right) - Re\left( \frac{z{w}'\left(z \right) }{w\left(z \right) }\right) \\
		&\geq s +\sum_{j=1}^{m-s}\frac{z}{z-b_{j}}- Re\left( \frac{z{w}'\left(z \right) }{w\left(z \right) }\right) .
	\end{align*}
	For $z \in D_{1}$, this gives with the help of Lemma \ref{4.2} that 
	\begin{align*}
		Re\left(\frac{z{r}'\left(z \right) }{r\left(z \right) } \right) 
		&\geq  s+ \sum_{j=1}^{m-s}\frac{z}{z-b_{j}}  - \left( \frac{n-\left|{B}'\left( z \right)  \right| }{2}\right) \\
		& \geq s+ \sum_{j=1}^{m-s}\frac{1}{1+\left|b_{j} \right| } - \left( \frac{n-\left|{B}'\left( z \right)  \right| }{2}\right) \\
			& \geq \sum_{j=1}^{m-s}\frac{1}{1+\left|b_{j} \right| } -  \frac{n}{2}+\frac{\left|{B}'\left( z \right)  \right| }{2} +\frac{m+sk}{1+k}-\frac{m-s}{1+k}\\
		& =\frac{1}{2}\left\lbrace \left|{B}'\left(z \right)  \right| -n+\frac{2(m+sk)}{1+k}+ 2\left(\sum_{j=1}^{m-s}\frac{1}{1+\left|b_{j} \right| }-\frac{m-s}{1+k} \right)\right\rbrace .
	\end{align*}
	This for $z \in D_{1},$ straightforwardly then
	\begin{align*}
		\left|\frac{{r}'\left(z \right) }{r\left( z\right) } \right|
		& = \left|\frac{z{r}'\left(z \right) }{r\left( z\right) } \right|\geq Re \left(\frac{z{r}'\left(z \right) }{r\left( z\right) }\right) \\
		&\geq  \frac{1}{2}\left\lbrace \left|{B}'\left(z \right)  \right| -n+\frac{2(m+sk)}{1+k}+ 2\left(\sum_{j=1}^{m-s}\frac{1}{1+\left|b_{j} \right| }-\frac{m-s}{1+k} \right)\right\rbrace.  	
	\end{align*}

	Which proves the Theorem \ref{thm3.1}.
\end{proof}

\begin{proof}[\textbf{Proof of Corollary 3.2}] 
	
	From inequality (\ref{3.1}), we have
	\begin{align}{\label{3.4}}
		\left|{r}'(z) \right| \geq\frac{1}{2}\left\lbrace \left|{B}'(z)\right|- n + \frac{2(m+sk)}{1+k} +2\left(\sum_{j=1}^{m-s}\frac{1}{1+|b_{j}|}-\frac{m-s}{1+k} \right)  \right\rbrace \left|r(z) \right| \nonumber\\
		= \frac{1}{2}\left\lbrace \left|{B}'(z)\right|- n + \frac{2(m+sk)}{1+k} +2\sum_{j=1}^{m-s}\left( \frac{1}{1+|b_{j}|}-\frac{1}{1+k} \right)  \right\rbrace \left|r(z) \right| 	\nonumber\\ 
		= \frac{1}{2}\left\lbrace \left|{B}'(z)\right|- n + \frac{2(m+sk)}{1+k} +\frac{2k}{k+1}\sum_{j=1}^{m-s}\frac{k-|b_{j}|}{k+k|b_{j}|}  \right\rbrace \left|r(z) \right| 	\nonumber\\ 
		\geq  \frac{1}{2}\left\lbrace \left|{B}'(z)\right|- n + \frac{2(m+sk)}{1+k}+ \frac{2k}{k+1}\sum_{j=1}^{m-s}\frac{k-|b_{j}|}{k+|b_{j}|}  \right\rbrace \left|r(z) \right| 	\nonumber\\
		= \frac{1}{2}\left\lbrace \left|{B}'(z)\right|- n + \frac{2(m+sk)}{1+k}+ \frac{2k}{k+1}\sum_{j=1}^{m-s}\frac{1-\frac{|b_{j}|}{k}}{1+\frac{|b_{j}|}{k}}  \right\rbrace \left|r(z) \right| 	
	\end{align}
	Since $\frac{|b_{j}|}{k} \leq 1 ,$ by Lemma \ref{4.4} , we conclude from inequality (\ref{3.4}) that
	\begin{align*}
		\left|{r}'(z) \right| \geq  \frac{1}{2}\left\lbrace \left|{B}'(z)\right|- n + \frac{2(m+sk)}{1+k}+ \frac{2k}{k+1}\left( \frac{1-\sum_{j=1}^{m-s}\frac{|b_{j}|}{k}}{1+\sum_{j=1}^{m-s}\frac{|b_{j}|}{k}}\right)   \right\rbrace \left|r(z) \right| \\
		= \frac{1}{2}\left\lbrace \left|{B}'(z)\right|- n + \frac{2(m+sk)}{1+k}+ \frac{2k}{k+1}\left( \frac{k^{m}|c_{m}|-|c_{s}|}{k^{m}|c_{m}|+|c_{s}|}\right)   \right\rbrace \left|r(z) \right|.
	\end{align*}
\end{proof}

	\begin{proof}[\textbf{Proof of Theorem 3.2}] Assume  that $r\in \Re_{n}$ has no zeros in $D_{k}^{+}$ where $k\leq 1$.
	Let $m^{'}=\min_{\left| z\right|=k }\left|r(z) \right|,$ then 
	$m^{'} \leq \left|r\left(z \right)  \right| $ for $z \in D_{k}$.
	If $r\left( z\right) $ has a zero on $\left|z \right|=k, $ then $m^{'}=0,$ hence for every $\alpha$ with $\left|\alpha \right| < 1$ we get $r\left(z \right)+\alpha m^{'} =r\left( z \right). $
	In case $r\left(z \right) $ has no zeros on $\left| z\right|=k, $ we have for every $\alpha$ with  $\left|\alpha \right| < 1$ that 
	$\left|\alpha m  \right|<\left|r\left(z \right)  \right|  $
	for $\left| z\right|=k. $	
	It follows by Rouch\'{e}'s theorem that
	$R\left(z \right)=r\left(z \right) + \alpha m^{'}  $ and $r\left( z\right) $ have same number of zeros in $D_{k}^{-},$ that is, for every $\alpha$ with $\left|\alpha \right| < 1,$
	$R\left( z\right) $ has no zero in $D_{k}^{+}.$
	If $b_{1},b_{2},\dots,b_{m}$ are zeros of $R\left( z\right),$ $ m\leq n $
	and $\left|b_{j} \right|\leq k\leq 1, $ we have 
	\begin{align*}
		\frac{z{R}'\left(z \right) }{R\left(z \right) }
		=  \frac{z{r}'\left(z \right) }{r\left(z \right) }
		&=	 \frac{z{p}'\left(z \right) }{p\left(z \right) } - \frac{z{w}'\left(z \right) }{w\left(z \right) }\\
		&=\sum_{j=1}^{m}\frac{z}{z-b_{j}}- \frac{z{w}'\left(z \right) }{w\left(z \right) }.
	\end{align*}
	For $z \in D_{1}$, gives with the help of Lemma \ref{4.2}, that 
	\begin{align*}
		Re\left(\frac{z{R}'\left(z \right) }{R\left(z \right) } \right) 
		&= Re \left(\sum_{j=1}^{m}\frac{z}{z-b_{j}} \right)  - \left( \frac{n-\left|{B}'\left( z \right)  \right| }{2}\right) \\
		& \geq \sum_{j=1}^{m}\frac{1}{1+\left|b_{j} \right| } - \left( \frac{n-\left|{B}'\left( z \right)  \right| }{2}\right) \\
		& = \frac{\left|{B}'\left(z \right)  \right| }{2} + \left(\sum_{j=1}^{m}\frac{1}{1+\left|b_{j} \right| }-\frac{n}{2} \right)+\frac{n-m}{1+k}-\frac{n-m}{1+k}\\
			& = \frac{\left|{B}'\left(z \right)  \right| }{2} + \frac{n(1-k)-2m}{2(1+k)}+ \left(\sum_{j=1}^{m}\frac{1}{1+\left|b_{j} \right| }-\frac{n-m}{1+k} \right),
	\end{align*}
	where $R\left(z \right)\neq 0. $ Then
	\begin{align*}
		\left|\frac{{R}'\left(z \right) }{R\left( z\right) } \right|
		& = \left|\frac{z{R}'\left(z \right) }{R\left( z\right) } \right|\geq Re \left(\frac{z{R}'\left(z \right) }{R\left( z\right) }\right)\\ 
		& \geq  \frac{\left|{B}'\left(z \right)  \right| }{2} + \frac{n(1-k)-2m}{2(1+k)}+ \left(\sum_{j=1}^{m}\frac{1}{1+\left|b_{j} \right| }-\frac{n-m}{1+k} \right).  	
	\end{align*}
	This implies that
	$$ \left|{R}'\left(z \right)  \right|\geq \left[ \frac{\left|{B}'\left(z \right)  \right| }{2} + \frac{n(1-k)-2m}{2(1+k)}+ \left(\sum_{j=1}^{m}\frac{1}{1+\left|b_{j} \right| }-\frac{n-m}{1+k} \right) \right]\left|R\left(z \right)  \right|, \  \text{for}\  \ z \in D_{1}.$$
	Since $R\left(z \right)=r\left( z\right)+ \alpha m^{'} , $ we get
	$$ \left|{r}'\left(z \right)  \right|\geq \left[ \frac{\left|{B}'\left(z \right)  \right| }{2} + \frac{n(1-k)-2m}{2(1+k)}+ \left(\sum_{j=1}^{m}\frac{1}{1+\left|b_{j} \right| }-\frac{n-m}{1+k} \right) \right]\left|r\left(z \right) +\alpha m^{'} \right|,\
	  \text{for} \ z\in D_{1}.$$
	Note that this inequality is trivially true for $R\left(z \right)=0.$ Therefore, this inequality holds for all $z \in D_{1}.$ Choosing the argument of $\alpha$ suitably in the right side of the above inequality and noting that the left side is independent of $\alpha$, we get that 
	$$ \left|{r}'\left(z \right)  \right|\geq \left[ \frac{\left|{B}'\left(z \right)  \right| }{2} + \frac{n(1-k)-2m}{2(1+k)}+ \left(\sum_{j=1}^{m}\frac{1}{1+\left|b_{j} \right| }-\frac{n-m}{1+k} \right)  \right]\left( \left|r\left(z \right)\right|  +\left| \alpha\right|  m^{'}  \right),$$ for $z\in D_{1}.$
	Letting $\left|\alpha \right|\rightarrow 1 $, we get for $z \in D_{1}$, that
	\begin{align*}
		\left|{r}'\left(z \right)  \right|
		&\geq \frac{1}{2}\left[ \left|{B}'\left(z \right)  \right|  + \frac{n(1-k)-2m}{1+k}+ 2\left(\sum_{j=1}^{m}\frac{1}{1+\left|b_{j} \right| }-\frac{n-m}{1+k} \right)  \right]\\
		&\times\left( \left|r\left(z \right)\right|  + m^{'} \right)\\
		&= \frac{1}{2}\left[ \left|{B}'\left(z \right)  \right|  + \frac{n(1-k)-2m}{1+k}+ 2\left(\sum_{j=1}^{m}\frac{1}{1+\left|b_{j} \right| }-\frac{n-m}{1+k} \right)    \right]\\
		&\times\left( \left|r\left(z \right)\right|  + m^{'} \right),
	\end{align*}
	which proves the Theorem \ref{thm3.2}.
\end{proof}
\textbf {Declaration :}\\

\textbf {Ethical Approval :} 
	Not applicable for the present work. It deals neither with human and/ nor animals. \\
	
\textbf {Competing interests :} 
	Authors have no competing interests to disclose.\\

\textbf {Funding :}
	This work received no specific grant from any funding agency in the public, commercial, or non profit sectors. \\
	
\textbf {Availability of data and materials :} 
	Data sharing is not applicable to this article as no data sets are generated or analysed during the current study.\\

\end{document}